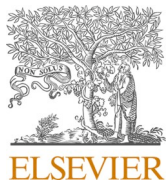
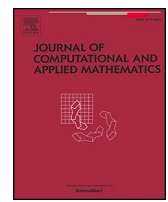

# An efficient algorithm for the minimal least squares solution of linear systems with indefinite symmetric matrices

Ibai Coria [a,*], Gorka Urkullu [a], Haritz Uriarte [b], Igor Fernández de Bustos [b]

[a] *Department of Applied Mathematics, Universidad del País Vasco/Euskal Herriko Unibertsitatea UPV/EHU, Bilbao 48013, Spain*
[b] *Department of Mechanical Engineering, Universidad del País Vasco/Euskal Herriko Unibertsitatea UPV/EHU, Bilbao 48013, Spain*

ABSTRACT

In this work, a new algorithm for solving symmetric indefinite systems of linear equations is presented. It factorizes the matrix into the form $LDL^t$ using Jacobi rotations in order to increase the pivot́s absolute value. Furthermore, Rooḱs pivoting strategy is also adapted and implemented. In determinate compatible systems, the computational cost of the algorithm was similar to the cost of the Bunch-Kaufman method, but the error was approximately 50 % smaller for intermediate and large matrices, regardless of the condition number of the coefficient matrix. Furthermore, unlike Bunch-Kaufman, the new algorithm calculates with little additional cost the fundamental basis of the null space, and obtains the minimal least squares and minimum norm solutions. In minimal least squares with minimum norm problems, the new algorithm was compared with the LAPACK Complete Orthogonal Decomposition algorithm, among others. The obtained error with both algorithms was similar but the computational cost was at least 20 % smaller with the new algorithm, even though the Complete Orthogonal Decomposition is implemented in a blocked form.

## 1. Introduction

Linear systems with symmetric coefficient matrices are commonly solved using Choleskys $LL^t$ factorization due to its high efficiency. However, the use of this factorization is limited because it is only valid for regular positive definite matrices. $LDL^t$ factorization is valid for both, positive and negative definite matrices, but not for the general case of indefinite matrices. For symmetric indefinite systems of linear equations, Gaussian elimination can be used, but it does not take advantage of the symmetry to reduce cost and storage. Therefore, for these cases, the Bunch-Kaufman factorization (also known as the diagonal pivoting method with partial pivoting) is generally used [1]. In fact, this algorithm is employed in both LINPACK and LAPACK. This algorithm factors $A$ as $A = PLDL^tP^t$, but in this case, $D$ is a block diagonal composed of $1 \times 1$ and $2 \times 2$ blocks [2,3].

Bunch-Kaufman algorithm can be considered as a partial pivoting method, while the Bunch-Parlett algorithm follows the same procedure but using complete (or total) pivoting [4]. Obviously, with Bunch-Kaufman the cost is reduced because the number of comparisons is O($n^2$) instead of O($n^3$), but at the expense of a larger growth factor, and therefore, a lower stability.

If these methods are used to solve a linear systems of equations, the coefficient matrix must be regular, as they have been only developed to solve determinate compatible systems. Therefore, they do not provide minimal norm and minimal least squares solutions. For that case, the eigenvalues and eigenvectors of the matrix should be obtained in order to perform a spectral decomposition in the

---

form $A = UDU^t$, where $U$ is an orthogonal matrix [5]. As it is known, spectral decompositions provide very accurate results but at the expense of a higher cost [6].

In this work, an alternative to Bunch-Kaufman and spectral decomposition is developed, which uses Jacobi rotations in order to maximize the pivots, and thus, minimize the growth factor [5]. The new algorithm, which is based on a previous work by the authors [6,7], factorize the matrix, calculates its null space with little additional cost, and solves efficiently the minimal least squares problem known as $min_x \| b - Ax_2 \|$, with $\|x_2\|$ minimal, using a procedure that combines Peters-Wilkinson and Sautter methods [8,9]. Furthermore, the Rook's pivoting method used in Gaussian Elimination is adapted and implemented to further improve the relation between cost and accuracy. In Gaussian Elimination, this pivoting strategy searches the pivot that is maximal in absolute value in both its row (to minimize the growth of error in the sweepout phase) and column (to foster stability during the back-substitution phase). In fact, while partial pivoting might have exponential error growth in the worst-case scenario, it has been demonstrated that Rook pivoting strategy never has exponential error growth [10]. Thus, comparing to partial pivoting, it is a more accurate strategy with very little additional cost [11,12].

The work is organized as follows: In Sections 2 and 3 the algorithm of the factorization is described and Rook's pivoting strategy is implemented. In Section 4 the growth factor of the elements along the factorization is calculated, and in Section 5 the fundamental basis of the null space is easily calculated. Section 6 provides the minimal least squares solution with minimum norm, while Section 7 presents numerical experiments. Finally, some conclusions are drawn.

## 2. Algorithm description

The presented algorithm, in addition to rows and columns permutations, uses Jacobi rotations. The aim is to preserve the symmetry along the factorization process in order to keep the computational and storage cost as low as possible. Jacobi rotations have been generally used to factor the matrix in orthogonal and diagonal matrices, but not for matrix triangularizations [13]. For simplicity, if row/column permutations are omitted, Jacobi rotations modify the first two rows and columns in such a way that zeroes element $a_{12}$ and increase the absolute value of the pivot $a_{11}$. For that purpose, a rotation of angle $\theta$ is applied. The calculation of this angle is clearly explained in [5,7]. Therefore, taking into account the nomenclature used in Eq. (1), $|b_{11}| \geq |a_{11}|$, and therefore, the growth factor is minimized. After the rotation, element $b_{11}$ is used to perform the first step of the $LDL^t$ factorization.

$$\begin{pmatrix} \cos(\theta) & -\sin(\theta) & 0^t_{n-2} \\ \sin(\theta) & \cos(\theta) & 0^t_{n-2} \\ 0_{n-2} & 0_{n-2} & I_{n-2} \end{pmatrix} \begin{pmatrix} a_{11} & a_{12} & a^t_1 \\ a_{12} & a_{22} & a^t_2 \\ a_1 & a_2 & A_{(n-2)x(n-2)} \end{pmatrix} \begin{pmatrix} \cos(\theta) & \sin(\theta) & 0^t_{n-2} \\ -\sin(\theta) & \cos(\theta) & 0^t_{n-2} \\ 0_{n-2} & 0_{n-2} & I_{n-2} \end{pmatrix} = \begin{pmatrix} b_{11} & 0 & b^t_1 \\ 0 & b_{22} & b^t_2 \\ b_1 & b_2 & A_{(n-2)x(n-2)} \end{pmatrix}. \quad (1)$$

In this sense, as the first two rows and columns will have influence on the pivot value, first two permutations are carried out using a matrix $P$, and later a matrix $G$ is used to apply the explained rotation. Thus, a given matrix $A$ can be factored in the form

$$A = P_1 G_1 P_2 G_2 \ldots P_k G_k \ldots P_r G_r L D L^t G_r^t P_r^t \ldots G_k^t P_k^t \ldots G_2^t P_2^t G_1^t P_1^t, \quad (2)$$

where $r$ is the rank of the matrix and the index $k$ the number of an intermediate step. Note that in every step, the permutations are first applied, and next the rotation of the rows $k$ and $k+1$. In order to simplify Eq. (2), every permutation and rotation can be included in a single $M$ orthogonal matrix:

$$M = \prod_{k=1}^{r} P_k G_k. \quad (3)$$

Therefore, matrix $A$ can be factored as

$$A = MLDL^t M^t. \quad (4)$$

In contrast to the algorithm for symmetric indefinite matrices previously developed by the authors in [7], the permutations and rotations are applied to the whole matrix, and not only to the submatrix of the Schur complement. This considerably simplifies the calculation of the minimal least squares solution, because it makes possible to lay aside the permutation and rotation matrices during the null space calculation. Thus, the structure of the matrices is maintained.

Row major order in which $j \geq i$ is used, a $2n$ sized integer vector has to be arranged to store the pivot indexes, and for the rotations a $n$ sized double vector is used. Note that, instead of the rotated angle, the tangent of the rotated angle is stored to avoid the use of trigonometric functions [5,7].

## 3. Modified Rook's pivoting method

As explained, Rook's method is a pivoting strategy that is mainly used in Gaussian Elimination which searches the largest pivot taking into account both its row and column. In the presented algorithm, the rotation is carried out with the first two rows and columns, and as a consequence of the rotation, a large element in the second row could generate a large element in the first row. Therefore, in order to avoid those large elements, the pivot of step $k$ should be the largest in $k$ and $k+1$ rows (note that columns are equal due to the symmetry of the system).

Let $a_{ij}$ be the element selected with the pivoting strategy in the step $k$. If the pivot $i \neq j$, the absolute value of elements $a_{ii}$ and $a_{jj}$ have to be compared, in order to place the largest element in $(k,k)$ position. Thus, if $|a_{ii}| \geq |a_{jj}|$, $piv[2k] = i$ and $piv[2k+1] = j$; else, $piv[2k] = j$ and $piv[2k+1] = i$. With those permutations, $a_{ij}$, which is the element selected with the pivoting strategy (and therefore, the largest element in rows $i$ and $j$), will be placed in $(k, k+1)$ position. Afterwards, with the rotation, it will become zero and it will increase even more the absolute value of the element positioned in $(k,k)$ position.

If $i = j$, $piv[2k] = i$ and thus, $a_{ii}$ its moved to the $(k,k)$ position. In this case, another row with every element smaller than $|a_{ii}|$ is necessary in order to perform the second permutation. In order to minimize cost, the row studied in the previous iteration is chosen. Thus, according to the Algorithm 1, $piv[2k+1] = previous$. For the sake of clarity, but without loss of generality, in Algorithm 1 it is assumed that the whole symmetric matrix is stored, so the storage method of the symmetric matrix has been omitted. The algorithm returns the pivot value, its row and column, and for the case $i = j$, the previously studied row (*max*, *rmax*, *cmax* and *previous*, respectively).

## 4. Bound on the growth of elements

Large growth factors are undesirable because they warn about a likely numerical instability in the factorization, as shown by Wilkinson [14,15]. In this sense, the aim of this section is to bound the maximum growth factor that may suffer the matrix along the factorization. The growth factor can be defined as

**Algorithm 1**
Modified Rook's pivoting strategy.

```
max=0
for i = 1: n {full pivoting is applied until a non zero element is found}
  for j = i: n
    if (|a_ij| > max)
      continue = 1 {This is the stop criterion of the while loop}
      max, rmax, cmax ← |a_ij|, i, j
      if (i = j) {the pivot is in the diagonal, so the previous row is saved}
        continue = 0
        previous = i - 1
      end if
    end if
  end for
  if max > tolerance
    break
  end if
end for
if (max < tolerance) {if this condition is fulfilled, the submatrix being studied is null}
  return
end if
if (rmax = cmax & rmax = 1) {if a_11 is the largest element of the first row, the second row must be also studied}
  continue = 0
  previous = 2
  for j = 2: n
    if (|a_2j| > max)
      continue = 1
      max, rmax, cmax ← |a_2j|, 2, j
  if (2 = j)
        continue = 0
        previous = 1
      end if
    end if
  end for
end if
while (continue = 1) {Loop until an element is found which is the largest in its row and column}
continue = 0
  previous = i
  i = cmax
  for j = 1: n
    if (|a_ij| > max)
      continue = 1
      max, rmax, cmax ← |a_ij|, i, j
      if (i = j)
        continue=0
      end if
    end if
  end for
end while
```

$$\rho = \frac{max_{i,j,k}\left|a_{i,j}^{(k)}\right|}{max_{i,j}\left|a_{i,j}\right|}, \tag{5}$$

where $a_{i,j}^{(k)}$ are the elements of the matrix at the $k^{th}$ step. To calculate the bound of the growth factor, first the reduced matrices of the Schur complement will be defined as $A^{(k)}$, which can be factored as $A^{(k)} = L^{(k)}D^{(k)}L^{(k)^t}$. $L^{(k)}$ is a lower triangular matrix, and even Rooks pivoting strategy is used, due to the Jacobi rotations $max_{i,j}\left|L_{i,j}^{(k)}\right| \leq \sqrt{2}$. Note that in Gaussian elimination with Rook pivoting strategy or in Cholesky factorization, the elements value would be equal or smaller than 1. Thus, taking into account that the pivot can be $\sqrt{2}$ times smaller than the elements of its row, and that the elements of matrix $A$ are normalized so that $max_{i,j}|a_{i,j}| = 1$, the maximum value of the elements of $A^{(k)}$ can be defined as

$$max_{i,j}\left|a_{i,j}^{(k)}\right| \leq 1 + 2\sum_{i=1}^{k-1} d_i, \tag{6}$$

where $d_i$ are the elements of $D$ once the matrix is factored, that is to say, the pivots. Using Hadamard's inequality the next expression is obtained:

$$\left|Det\left(A^{(k)}\right)\right| \leq k^{k/2} \cdot \left(max_{i,j}\left|a_{i,j}^{(k)}\right|\right)^k \leq k^{k/2} \cdot \left(1 + 2\sum_{i=1}^{k-1} d_i\right). \tag{7}$$

Next, the elements of $D$ are used to calculate the determinant of $A^{(k)}$. Without loss of generality, $d_n \geq d_{n-1} \geq \ldots \geq d_1 \geq 0$ will be considered:

$$Det\left(A^{(k)}\right) = d_n d_{n-1} \ldots d_{n+1-k}. \tag{8}$$

Therefore:

$$d_n d_{n-1} \ldots d_{n+1-k} \leq k^{k/2} \cdot \left(1 + 2\sum_{i=1}^{k-1} d_i\right) \text{ for } k = 1, 2, \ldots, n. \tag{9}$$

The maximum growth factor is achieved when all inequalities in (9) are equalities, and in that case, as the elements of matrix $A$ are normalized, $\rho = d_n$. Following a procedure similar to the proof of Theorem 1 presented in [10], the next equation can be deducted to calculate the growth factor:

$$\rho \leq s_1(1 + 2s_2)(1 + 2s_3) \cdots (1 + 2s_n), \tag{10}$$

where

$$s_k(1 + 2s_k)^{k-1} = \frac{k^{k/2}}{(k-1)^{(k-1)/2}}. \tag{11}$$

Nevertheless, Eq. (11) must be solved numerically, which can be tedious because of the difficulty of comparing the growth factor of different factorization methods. Therefore, in order to obtain an analytical expression, a procedure similar to the proof of Theorem 3 presented in [10] can be performed. Thus, it can be proved that the next equation is always fulfilled:

$$\rho \leq 2.8 n^{3\ln(n)/4}. \tag{12}$$

This bound is less than the double of the Gaussian Elimination with the Rook's pivoting strategy bound, which can be considered a very good result. Furthermore, for large matrices, the grow factor is more slowly growing than any function of the form $ab^n$ for $a > 0$ and $b > 1$. Instead, Bunch-Kaufman factorization bound is $(2.57)^{n-1}$, which means that it might have exponential error growth [16].

## 5. Calculation of the null space

In this section, the null space of the permutated and rotated matrix will be calculated, since it will be necessary to obtain the minimal least squares solution. If the matrices are divided in blocs according to the rank $r$, and the matrix dimension is $n$:

$$LDL^tN = \begin{pmatrix} L_{11} & (0) \\ L_{21} & I_{n-r} \end{pmatrix} \begin{pmatrix} D_{11} & (0) \\ (0) & (0) \end{pmatrix} \begin{pmatrix} L_{11}^t & L_{21}^t \\ (0) & I_{n-r} \end{pmatrix} \begin{pmatrix} N_1 \\ N_2 \end{pmatrix} = (0). \tag{13}$$

As $L$ is a regular matrix:

$$\begin{pmatrix} D_{11} & (0) \\ (0) & (0) \end{pmatrix} \begin{pmatrix} L_{11}^t & L_{21}^t \\ (0) & I_{n-r} \end{pmatrix} \begin{pmatrix} N_1 \\ N_2 \end{pmatrix} = (0). \tag{14}$$

The second row is null, so it does not provide any information. Developing the first row:

$$D_{11}L_{11}{}'N_1 + D_{11}L_{21}{}'N_2 = (0), \tag{15}$$

$$N_1 = -(L_{11}{}')^{-1}L_{21}{}'N_2. \tag{16}$$

Every matrix $N_2$ will provide a valid expression of the null space. Therefore, the next basis of the null space can be used:

$$N = \begin{pmatrix} -(L_{11}{}')^{-1}L_{21}{}^t \\ I_{n-r} \end{pmatrix}. \tag{17}$$

The obtained null space expression corresponds to the fundamental null basis explained in [17]. The fundamental basis has three main advantages: the computational cost required for its calculation is very competitive since $L_{11}$ is a triangular matrix, it is not necessary to store the last $n - r$ rows since it will always be the identity matrix, and multiplying it with another matrix entails a lower number of operations. Furthermore, in order to reduce storage, $(L_{11}{}^t)^{-1}L_{21}{}^t$ can be stored in the place of $L_{21}{}^t$ (which is usually the place where $A_{21}$ was stored). Note that, once the null space is calculated, $L_{21}{}^t$ will no longer be necessary.

## 6. Minimal least squares solution with minimum norm

To obtain the minimal least squares solution with minimum norm the procedure used by the authors in [6] will be employed. Let's assume a symmetric linear system of equations with $r \leq n$:

$$Ax = b. \tag{18}$$

To obtain the minimal least squares solution $\hat{x}$, the normal equations are used:

$$AA\hat{x} = Ab, \tag{19}$$

$$(MLDL^T M^T)(MLDL^T M^T)\hat{x} = (MLDL^T M^T)b, \tag{20}$$

$$DL^T LDL^T M^T \hat{x} = DL^T M^T b. \tag{21}$$

This system is, in the more general case, an indeterminate system. In order to solve the system, two changes have to be introduced:

$$z = DL^T M^t \hat{x}, \tag{22}$$

$$c = M^t b. \tag{23}$$

So one needs to solve

$$DL^T Lz = DL^T c. \tag{24}$$

With the partition used in Section 5:

$$\begin{pmatrix} D_{11} & (0) \\ (0) & (0) \end{pmatrix} \begin{pmatrix} L_{11}{}^t & L_{21}{}^t \\ (0) & I_{n-r} \end{pmatrix} \begin{pmatrix} L_{11} & (0) \\ L_{21} & I_{n-r} \end{pmatrix} \begin{pmatrix} z_1 \\ z_2 \end{pmatrix} = \begin{pmatrix} D_{11} & (0) \\ (0) & (0) \end{pmatrix} \begin{pmatrix} L_{11}{}^t & L_{21}{}^t \\ (0) & I_{n-r} \end{pmatrix} \begin{pmatrix} c_1 \\ c_2 \end{pmatrix}. \tag{25}$$

The last $n - r$ equations are null, and from the first $r$ equations the next system is obtained:

$$(D_{11}L_{11}{}^T L_{11} + D_{11}L_{21}{}^T L_{21})z_1 = D_{11}L_{11}{}^T c_1 + D_{11}L_{21}{}^T c_2 - D_{11}L_{21}{}^T z_2. \tag{26}$$

As any solution of the system is valid:

$$(L_{11}{}^T L_{11} + L_{21}{}^T L_{21})z_1 = L_{11}{}^T c_1 + L_{21}{}^T c_2, \tag{27}$$

$$z_2 = 0. \tag{28}$$

Finally, in order to take advantage of the null space, the previous system can be expressed as

$$(I + N_1 N_1{}^T)L_{11}z_1 = c_1 - N_1 c_2. \tag{29}$$

The matrix $(I + N_1 N_1{}^T)$ of the system, whose dimension is $r$, is positive definite and is well conditioned. Therefore, to obtain $L_{11}z_1$, PCG or Cholesky methods can be used. As explained in [6], this approach is quite useful when dealing with matrices of low rank, because the obtained system will be quite small. If $r > (n/2)$, a different approach must be used for the calculation of $z_1$. Let's assume the change

$$A\hat{x} = Mu. \tag{30}$$

Going back to Eq. (20):

$$MLDL^T M^T Mu = MLDL^T M^T b, \tag{31}$$

$$DL^T u = DL^T M^T b. \tag{32}$$

Applying the change of Eq. (23):

$$DL^T u = DL^T c. \tag{33}$$

Using the previously calculated null space, the general solution of the system can be obtained:

$$u = u_p + N\alpha = c + N\alpha. \tag{34}$$

The obtained solution has to satisfy also Eq. (30), so $N^T u$ must be null:

$$N^t u = 0 = N^t c + N^t N\alpha, \tag{35}$$

$$N^t N\alpha = -N^t c, \tag{36}$$

$$(N_1^t N_1 + I_{n-r})\alpha = -N_1^t c_1 - c_2, \tag{37}$$

which has again a positive definite and well conditioned matrix, but, in this case, matrix dimension is $(n-r)x(n-r)$. Once $\alpha$ is known, from Eq. (34) the first $r$ elements of vector $u$ are calculated:

$$u_1 = c_1 + N_1 \alpha. \tag{38}$$

Finally, going back to Eq. (30):

$$LDL^t M^t \widehat{x} = u, \tag{39}$$

$$DL^t M^t \widehat{x} = z = L^{-1} u. \tag{40}$$

Once again, $z_2$ can be assumed as the null vector, and thus, $z_1$ is obtained as

$$z_1 = L_{11}^{-1} u_1. \tag{41}$$

Once $z$ is obtained, the minimum norm solution $\widehat{x}_{min}$ must be obtained. For that purpose, first the next change will be introduced in Eq. (40):

$$w = M^t \widehat{x}, \tag{42}$$

$$DL^t w = z, \tag{43}$$

$$\begin{pmatrix} D_{11} & (0) \\ (0) & (0) \end{pmatrix} \begin{pmatrix} L_{11}^t & L_{21}^t \\ (0) & I_{n-r} \end{pmatrix} \begin{pmatrix} w_1 \\ w_2 \end{pmatrix} = \begin{pmatrix} z_1 \\ z_2 \end{pmatrix}. \tag{44}$$

Again, the last $n-r$ equations are null. Thus, it is possible to write the equation in the form

$$L_1 = \begin{pmatrix} L_{11} \\ L_{21} \end{pmatrix}, \tag{45}$$

$$L_1^t w = D_{11}^{-1} z_1. \tag{46}$$

To solve the system, the usual minimum norm expression has to be used:

$$w_{min} = \begin{pmatrix} w_{min1} \\ w_{min2} \end{pmatrix} = L_1 (L_1^t L_1)^{-1} D_{11}^{-1} z_1. \tag{47}$$

Substituting and simplifying the equation:

$$(I + N_1 N_1^t) w_{min1} = (L_{11}^t)^{-1} D_{11}^{-1} z_1, \tag{48}$$

$$w_{min2} = -N_1^T w_{min1}. \tag{49}$$

The system of Eq. (48) has a dimension $r$, and as it is again positive definite, it can be solved using Cholesky or PCG. As in the calculation of $z$, if $r > (n/2)$, a different approach must be used for the calculation of $w_{min}$. The general solution of Eq. (46) can be also

expressed as

$$w = w_p + N\beta, \tag{50}$$

$$w_p = \begin{pmatrix} w_{p1} \\ w_{p2} \end{pmatrix} = \left\{ \begin{array}{c} (L_{11}{}^T)^{-1} D_1{}^{-1} z_1 \\ 0 \end{array} \right\}. \tag{51}$$

Introducing the minimal norm condition:

$$N^t N \beta = -N^t w_p, \tag{52}$$

$$(N_1{}^T N_1 + I)\beta = -N_1{}^T w_{p1}. \tag{53}$$

Once $\beta$ is obtained, from Eq. (50) $w_{min}$ is calculated as

$$w_{min1} = w_{p1} + N_1 \beta, \tag{54}$$

$$w_{min2} = \beta. \tag{55}$$

Finally, to conclude with the algorithm, $\widehat{x}_{min}$ is calculated from Eq. (42):

$$\widehat{x}_{min} = M w_{min}. \tag{56}$$

In order to reduce storage as much as possible, the algorithm was implemented following the next steps (Algorithm 2).

## 7. Numerical results

The aim of the experiments was to solve linear systems of equations with the form $Ax = b$. They were carried out in a E5–2637 v3 @3.5 GHz (Haswell architecture), forcing the CPU to run at 3.5 GHz. A single processor was used but obviously including SIMD extension, in this case AVX2. The OS is Ubuntu Linux 22, and openblas implementation for the BLAS and LAPACK routines was employed, recompiling the libraries with gcc, using options "-march = haswell – O3" to take the most out of the processor. In order to do a proper comparison, single threaded openblas was used.

First, Bunch-Kaufman, Bunch-Parlett and algorithm developed by the authors in [7] (which uses in a similar way rotation matrices

---

**Algorithm 2**
Minimal least squares solution with minimum norm.

**Step 1.** Perform $A = MLDL^t M^t$ factorization using Rook's pivoting strategy.
$A \leftarrow D$ and $L^t$
$p, r \leftarrow M$
Where $p$ and $r$ are the permutation and rotation vectors, with size $2n$ (integer) and $n$ (double), respectively.
**Step 2.** Compute the null space.
$A_{12} \leftarrow (L_{11}{}^t)^{-1} L_{21}{}^t$
**Step 3.** Solve the minimal least squares problem.
$b \leftarrow c = M^t b$
if $r \leq (n/2)$:
  $b_1 \leftarrow c_1 - N_1 c_2 = b_1 + A_{12} b_2$
  $b_1 \leftarrow L_{11} z_1 = (I + N_1 N_1{}^T)^{-1}(c_1 - N_1 c_2) = (I + A_{12} A_{12}{}^T)^{-1} b_1$
else:
  $b_2 \leftarrow -N_1{}^t c_1 - c_2 = A_{12}{}^t b_1 - b_2$
  $b_2 \leftarrow \alpha = (N_1{}^t N_1 + I_{n-r})^{-1}(-N_1{}^t c_1 - c_2) = (A_{12}{}^t A_{12} + I_{n-r})^{-1} b_2$
  $b_1 \leftarrow u_1 = c_1 + N_1 \alpha = b_1 - A_{12} b_2$
end if
$b_1 \leftarrow z_1 = (L_{11})^{-1} b_1$
**Step 4.** Solve the minimum norm problem
$b_1 \leftarrow D_{11}{}^{-1} z_1 = D_{11}{}^{-1} b_1$
$b_1 \leftarrow (L_{11}{}^t)^{-1} D_{11}{}^{-1} z_1 = (L_{11}{}^t)^{-1} b_1$
if $r \leq (n/2)$:
  $b_1 \leftarrow w_{min1} = (I + N_1 N_1{}^t)^{-1}(L_{11}{}^t)^{-1} D_{11}{}^{-1} z_1 = (I + A_{12} A_{12}{}^t)^{-1} b_1$
  $b_2 \leftarrow w_{min2} = -N_1{}^T w_{min1} = A_{12}{}^t b_1$
else:
  $b_2 \leftarrow -N_1{}^T (L_{11}{}^t)^{-1} D_{11}{}^{-1} z_1 = A_{12}{}^t b_1$
  $b_2 \leftarrow w_{min2} = \beta = (N_1{}^T N_1 + I)^{-1}(-N_1{}^T (L_{11}{}^t)^{-1} D_{11}{}^{-1} z_1) = (A_{12}{}^T A_{12} + I)^{-1} b_2$
  $b_1 \leftarrow w_{min1} = (L_{11}{}^T)^{-1} D_1{}^{-1} z_1 + N_1 \beta = b_1 - A_{12} b_2$
end if
**Step 5.** Permute and rotate the elements to obtain the solution
$b \leftarrow \widehat{x}_{min} = M w_{min} = Mb$

**Table 1**
Reconstruction error in determinate compatible problems. Mean value (bold) and standard deviation (italics).

| Size | Bunch-Kaufman | Bunch-Parlett | Rotated CP | Rotated Rook |
| --- | --- | --- | --- | --- |
| 10 | **1.219E-15** | **9.476E-16** | **8.831E-16** | **1.098E-15** |
|  | *3.401E-16* | *2.282E-16* | *1.519E-16* | *1.879E-16* |
| 50 | **1.831E-14** | **1.163E-14** | **9.135E-15** | **1.158E-14** |
|  | *2.203E-15* | *7.681E-16* | *3.913E-16* | *5.731E-16* |
| 100 | **6.130E-14** | **3.626E-14** | **2.761E-14** | **3.517E-14** |
|  | *5.170E-15* | *1.380E-15* | *7.032E-16* | *1.105E-15* |
| 500 | **1.132E-12** | **5.897E-13** | **4.340E-13** | **5.695E-13** |
|  | *4.320E-14* | *6.950E-15* | *3.747E-15* | *7.272E-15* |
| 1 000 | **4.109E-12** | **2.070E-12** | **1.526E-12** | **2.018E-12** |
|  | *1.101E-13* | *1.521E-14* | *8.854E-15* | *1.725E-14* |

but in this case with full pivoting) were compared against the factorization algorithm developed in this work. In order to do so, symmetric regular matrices were used. The matrix $A$ and the solution vector $x$ were generated using the random number generator used in [7], with values $\mathbb{R}[-1,1]$. Next, $b = Ax$ was computed in quadruple precision. The results were averaged over 10 000 tests.

Table 1 shows the obtained mean (bold) and standard deviation (italics) of the reconstruction error with different factorization methods. The reconstruction error was calculated as the Frobenius norm of the original matrix minus the reconstructed matrix, this is, $\| A - MLDL^tM^t \|_F$. Quadruple precision was used to reconstruct the matrix in order to avoid rounding errors (all other calculations were carried out in double precision). The reconstruction error was used instead of the error of the linear system solution due to a lower scatter. Table 2 shows the approximate solution time in seconds needed for factorizing the matrix and solving the linear system. As can be appreciated, the accuracy of the new method (Rotated Rook) is similar to the Buch-Parlett method (which uses full pivoting) but with a similar cost to the Bunch-Kaufman method (which uses partial pivoting), thus showing the high efficiency of the new factorization method. Furthermore, the low scatter in both error and time shows the reliability of the obtained results.

In previous results, the influence of the condition number of matrix $A$ was not taken into account, so next, in order to test the algorithm in more severe conditions, the condition number of matrix $A$ was modified. For that purpose, the random matrices were generated through a spectral decomposition $A = UDU^t$. The maximum absolute value of the diagonal matrix was set as 1, and the minimum value was modified in order to obtain the desired condition number. The orthogonal matrix was built using the Stewart method [18] and the Mersenne Twister random generator along with the basic form of the Box-Mueller transformation [19–21]. Thus, random numbers with a Gaussian distribution in the range $\mathbb{R}[-1,1]$ were obtained. A random vector $z$ was also generated, which was assumed equal to $z = U^tb$. Therefore, $b = Uz$ and $x = UD^{-1}z$. Again, these initial calculations and the reconstruction error were computed in quadruple precision while all other calculations were carried out in double precision. Matrix size 100 was only studied due to the large number of parameters involved, and the results were averaged over 10 000 tests.

Table 3 shows the obtained reconstruction error and, as can be seen, the results are similar to the size 100 of Table 1. In this sense, it can be stated that the condition number does not have an important influence in the reconstruction error. Table 4 shows the error of the solution, which was calculated as the Euclidean norm of the exact solution minus the obtained solution, $\| x - x'_2 \|$. Obviously it increases as the condition number increases. However, the error of the new algorithm is always quite smaller than the one made by Bunch-Kaufman, and again is similar to the one made by Bunch-Parlett. Once again, this shows the good performance of the proposed algorithm. The computational cost was not included due to its similarity to the results of Table 2.

Finally, linear systems were solved using minimal least squares solutions with minimum norm. The random matrices were again generated as $A = UDU^t$, but in this case, matrix rank was controlled with the null elements of the diagonal of matrix $D$, while the non-zero elements of the diagonal were randomly generated with a Gaussian distribution within the range $\mathbb{R}[-1,1]$. The orthogonal matrix and the vector $z$ were generated following the previous procedure. Therefore, $b = Uz$ and $x = UD^+z$. Due to the fact that there are many parameters involved, to get a reasonable idea of the problem, it was assumed a rank equal to $n/2$, and $n/4$ incompatibilities in the problem.

**Table 2**
Approximate solution time (s) in determinate compatible problems. Mean value (bold) and standard deviation (italics).

| Size | Bunch-Kaufman | Bunch-Parlett | Rotated CP | Rotated Rook |
| --- | --- | --- | --- | --- |
| 10 | **3.219E-06** | **2.383E-06** | **2.605E-06** | **2.852E-06** |
|  | *1.072E-06* | *4.790E-07* | *5.321E-07* | *5.344E-07* |
| 50 | **2.759E-05** | **4.040E-05** | **4.091E-05** | **2.974E-05** |
|  | *1.751E-06* | *2.238E-06* | *1.654E-06* | *1.545E-06* |
| 100 | **1.105E-04** | **2.225E-04** | **2.331E-04** | **1.122E-04** |
|  | *7.784E-06* | *1.842E-05* | *9.800E-06* | *9.056E-06* |
| 500 | **8.425E-03** | **2.115E-02** | **2.462E-02** | **8.847E-03** |
|  | *6.010E-05* | *7.056E-04* | *5.398E-05* | *4.250E-05* |
| 1 000 | **6.435E-02** | **1.670E-01** | **1.985E-01** | **6.567E-02** |
|  | *2.804E-04* | *5.419E-03* | *4.398E-04* | *1.352E-04* |

**Table 3**
Reconstruction error in determinate compatible problems with different condition numbers. Mean value (bold) and standard deviation (italics).

| Condition number | Bunch-Kaufman | Bunch-Parlett | Rotated CP | Rotated Rook |
|---|---|---|---|---|
| 1E2 | **6.163E-15** | **3.525E-15** | **2.752E-15** | **3.471E-15** |
|  | *7.246E-16* | *2.911E-16* | *2.040E-16* | *2.627E-16* |
| 1E4 | **6.061E-15** | **3.469E-15** | **2.708E-15** | **3.418E-15** |
|  | *7.048E-16* | *2.910E-16* | *2.041E-16* | *2.603E-16* |
| 1E6 | **6.068E-15** | **3.464E-15** | **2.709E-15** | **3.419E-15** |
|  | *7.096E-16* | *2.958E-16* | *2.064E-16* | *2.643E-16* |
| 1E8 | **6.056E-15** | **3.465E-15** | **2.708E-15** | **3.418E-15** |
|  | *7.077E-16* | *2.922E-16* | *2.046E-16* | *2.612E-16* |
| 1E10 | **6.064E-15** | **3.463E-15** | **2.710E-15** | **3.423E-15** |
|  | *7.165E-16* | *2.894E-16* | *2.059E-16* | *2.616E-16* |

**Table 4**
Solution error in determinate compatible problems with different condition numbers. Mean value (bold) and standard deviation (italics).

| Condition number | Bunch-Kaufman | Bunch-Parlett | Rotated CP | Rotated Rook |
|---|---|---|---|---|
| 1E2 | **4.916E-25** | **2.401E-25** | **2.047E-25** | **2.258E-25** |
|  | *6.588E-25* | *3.076E-25* | *2.625E-25* | *2.919E-25* |
| 1E4 | **1.629E-17** | **8.192E-18** | **7.126E-18** | **7.561E-18** |
|  | *3.336E-17* | *1.716E-17* | *1.559E-17* | *1.644E-17* |
| 1E6 | **1.632E-09** | **8.103E-10** | **6.947E-10** | **7.520E-10** |
|  | *3.543E-09* | *1.700E-09* | *1.441E-09* | *1.610E-09* |
| 1E8 | **1.687E-01** | **7.916E-02** | **6.899E-02** | **7.651E-02** |
|  | *3.861E-01* | *1.659E-01* | *1.454E-01* | *1.602E-01* |
| 1E10 | **1.697E+07** | **7.922E+06** | **7.011E+06** | **7.786E+06** |
|  | *3.694E+07* | *1.670E+07* | *1.475E+07* | *1.683E+07* |

LAPACK does not have any algorithm to obtain minimal least squares solutions with minimum norm for symmetric matrices, so the new algorithm was compared with the complete SVD factorization (dgesvd), with the SVD factorization but saving time discarding the computation of the whole orthogonal matrices (dgelsd) and with the Complete Orthogonal Decomposition (dgelsy). It was also compared with the LDU algorithm developed by the authors in [6], which also uses Rook's pivoting strategy. The results were again averaged over 10 000 tests.

As can be appreciated in Table 5, for sizes above 52 the cost of the new method is much smaller than the previously developed LDU algorithm, which makes sense, since the new algorithm only works with the upper triangular part of the matrix. Regarding LAPACK algorithms, it can be seen that dgelsy is the most efficient; it has the lowest computational cost and a very competitive error. Nevertheless, if the new algorithm and dgelsy are compared, it can be appreciated that the new algorithm has a lower computational cost with also a very competitive error. It should be pointed out that the error was found to have a large dispersion, and therefore Table 6 just shows that the errors are of the same order of magnitude. The error of the solution was again calculated as the Euclidean norm of the exact solution minus the obtained solution, $\| x - x'_2 \|$.

Furthermore, two considerations must be taken into account regarding the computational cost. First, dgesvd, dgelsd and dgelsy are implemented with a blocked form, which delivers a lower computational cost in large matrices. Because of that, the new algorithm is more competitive for small matrices (when the whole matrix can be stored in the cache), than for large matrices, where a blocked form organizes better the data transfer from RAM memory to cache memory. In this sense, in a future work the presented algorithm could be also implemented in a blocked form in order to increase even more the efficiency of the algorithm. Second, the studied case where the rank is equal to $n/2$ and the problem has $n/4$ incompatibilities is the worst case for the presented algorithm, as explained in [6], due to the dimensions of the systems of Eqs. (29) and (48).

**Table 5**
Approximate solution time (s) in minimal least squares with minimum norm problems. Mean value (bold) and standard deviation (italics).

| Size | dgesvd | dgelsd | dgelsy | LDU Rook | Rotated Rook |
|---|---|---|---|---|---|
| 12 | **2.322E-05** | **2.095E-05** | **1.14E-05** | **4.660E-06** | **4.687E-06** |
|  | *3.416E-06* | *2.908E-06* | *1.194E-06* | *7.950E-07* | *7.172E-07* |
| 52 | **4.706E-04** | **2.710E-04** | **1.074E-04** | **4.935E-05** | **3.603E-05** |
|  | *1.963E-05* | *1.515E-05* | *4.654E-06* | *2.063E-06* | *3.905E-06* |
| 100 | **2.410E-03** | **1.007E-03** | **3.913E-04** | **2.548E-04** | **1.401E-04** |
|  | *6.606E-05* | *3.930E-05* | *6.859E-06* | *4.898E-06* | *1.052E-05* |
| 500 | **2.149E-01** | **4.236E-02** | **2.047E-02** | **3.447E-02** | **1.337E-02** |
|  | *3.449E-03* | *4.567E-04* | *1.277E-04* | *1.421E-04* | *1.823E-04* |
| 1000 | **1.824E+00** | **2.580E-01** | **1.294E-01** | **2.755E-01** | **1.070E-01** |
|  | *2.109E-02* | *1.967E-03* | *7.066E-04* | *4.756E-04* | *2.881E-04* |

**Table 6**
Error of the least squares solution with minimum norm. Mean value (bold) and standard deviation (italics).

| Size | dgesvd | dgelsd | dgelsy | LDU Rook | Rotated Rook |
| --- | --- | --- | --- | --- | --- |
| 12 | **4.594E-12** | **4.594E-12** | **8.023E-12** | **7.121E-12** | **6.066E-12** |
|  | *1.247E-10* | *1.247E-10* | *2.126E-10* | *2.493E-10* | *1.806E-10* |
| 52 | **3.454E-09** | **8.018E-08** | **6.337E-09** | **4.419E-09** | **2.829E-09** |
|  | *2.911E-07* | *7.773E-06* | *5.359E-07* | *3.244E-07* | *1.818E-07* |
| 100 | **5.280E-09** | **1.097E-07** | **1.397E-08** | **1.193E-08** | **1.399E-08** |
|  | *3.420E-07* | *6.593E-06* | *8.429E-07* | *8.012E-07* | *8.841E-07* |
| 500 | **1.420E-05** | **2.030E-05** | **1.422E-04** | **2.720E-04** | **1.330E-04** |
|  | *1.323E-03* | *1.873E-03* | *1.410E-02* | *2.698E-02* | *1.299E-02* |
| 1000 | **9.784E-07** | **1.004E-06** | **3.881E-06** | **1.083E-05** | **4.758E-06** |
|  | *6.850E-05* | *6.875E-05* | *2.667E-04* | *6.818E-04* | *2.295E-04* |

## 8. Conclusions

The presented algorithm factors symmetric indefinite matrices using Jacobi rotations along with a modified Rook's pivoting strategy. It can be considered an alternative to the commonly used Bunch-Kaufman method since, with a similar computational cost, it provides more accurate results due to the tight growth factor; in the performed tests, the error was approximately 50 % smaller for intermediate a large matrices, regardless the condition number of the coefficient matrix. The growth factor bound is less than the double of the Gaussian Elimination with the Rook's pivoting strategy bound, and therefore, the new algorithm will never have exponential error growth.

Furthermore, it has a very low memory cost, calculates the fundamental null basis with little additional cost, and obtains the minimal least squares and minimum norm solutions. In minimal least squares with minimum norm problems, the new algorithm and the Complete Orthogonal Decomposition algorithm had a similar error in the solution, but the computational cost of the new algorithm was at least 20 % smaller, even though the Complete Orthogonal Decomposition is implemented in a blocked form.

Therefore, numerical results show that the new procedure is more efficient than LAPACK algorithms in both compatible and incompatible systems, and even in ill-conditioned systems. Besides, the efficiency could be further improved in future works using blocked forms.

**Data availability**

Data will be made available on request.

**Acknowledgments**

The authors would like to thank to the Basque Government for its funding to the Research Group recognized under section IT1542-22. The authors also specially thank the grant PID2021-124677NB-I00 funded by MCIN/AEI/10.13039/501100011033 and by "ERDF A way of making Europe".